\documentclass[12pt,a4paper]{article}
\usepackage[cp1251]{inputenc}
\usepackage[english]{babel}
\usepackage[psamsfonts]{amssymb}
\begin{document}

%\large

\thispagestyle{empty}

\begin{center}

{\bf Twisted conjugacy classes of the unit element  \footnote[1]{ The authors were supported by “Development of the Scientific Potential of Higher School” (Project 2.1.1.10726),
the Interdisciplinary Project of the Siberian Division of the Russian Academy of Sciences (Grant No. 44–2012). and
the Federal Target Program “Scientific and Scientific-Pedagogical Personnel of Innovative Russia” for 2009–2013
(State Contract 02.740.11.5191).
\\}}

\bigskip

{V.~G.~Bardakov, T.~R. Nasybullov, M.~V.~Neshchadim}

\end{center}

\vspace{0.5cm}

\hfill{UDC 512.8}

\begin{quote}
\noindent{\bf Abstract. }
In this paper we study twisted conjugacy classes of the unit element in different groups.
A.~L.~Fel'shtyn and E.~V.~Troitsky showed that in an abelian group a twisted conjugacy class of the unit element is a subgroup for any automorphism. In this article we study the question about a structure of the groups, where twisted conjugacy class of the unit element is a subgroup for every automorphism (inner automorphism).


\end{quote}

\vspace{0.8cm}

\begin{center}

{\bf Introduction}

\end{center}

\vspace{0.5cm}

The conjugacy classes of the group reflect some properties of the group. For example, one of the questions in the combinatorial group theory is whether an infinite group must have an infinite number of conjugacy classes. In 1949 Higman-Neumann-Neumann \cite{HNN} constructed an infinitely generated group with the finite number of conjugacy classes. Later on, S.~Ivanov (\cite{Iv}, Theorem 41.2) constructed an example of the finitely generated group with this property.
Then D.~V.~Osin \cite{os} gave an example of the finitely generated infinite group where all non trivial elements are conjugated.

The generalization of the conjugacy classes are twisted conjugacy classes, namely $\varphi$-conjugacy classes, where $\varphi$ is an automorphism of the group.
Recently there have been a lot of papers devoted to the study of twisted conjugacy classes in different groups. The appropriate reviews can be found in \cite{fel, fel1, feltro2}.

The twisted conjugacy classes appear naturally in the Nielsen-Reidemeister fixed point theory.
Let $f : X \longrightarrow X$ be the map of the compact topological space $X$ on itslef,  $p : \widetilde{X} \longrightarrow X$ be the universal covering of $X$ and
$\widetilde{f} : \widetilde{X} \longrightarrow \widetilde{X}$ be the lifting of $f$, i.~e.
$p \circ \widetilde{f} = f \circ p$. Two liftings $\widetilde{f}$ and $\widetilde{f}'$ are cold {\it conjugated} if there is a such
$\gamma \in \pi_1(X)$ that $\widetilde{f}' = \gamma \circ \widetilde{f} \circ \gamma^{-1}$.
The subset $p(Fix(\widetilde{f}))$ is called {\it the fixed point class} of $f$, determined by the lifting class $[\widetilde{f}]$.
The fixed point class is called
{\it essential}, if its index is nonzero. The number of lifting classes of $f$ (and hence the
number of the fixed point classes) is called {\it the Reidemeister number} of $f$ and is denoted $R(f)$.
The number of the essential fixed point classes  is called {\it the Nielsen number} of $f$ and is denoted $N(f)$. $R(f)$ and $N(f)$ are homotopy type invariants.
The numbers $N(f)$ and $R(f)$ of the map $f$ are closely connected to each other and they are the main objects of investigation in the Nielsen-Reidemeister fixed point theory.

On the other hand, the map $f$ induces endomorphism $\varphi = f_\sharp$ of the fundamental group $\pi_1(X)$, where the number of $\varphi$-conjugacy classes of the endomorphism $\varphi$ denoted by symbol $R(\varphi)$, is equal to the Reqidemeister number $R(f)$, and is called {\it the Reidemeister number of the endomorphism} $\varphi$. Thus, the topological problem of calculating of the number $R(f)$ is reduced to a purely algebraic problem of calculatin of the number $R(\varphi)$ (more details in \cite{fel, fel1}).

As shown in the paper of A.~L.~Fel'shtyn and E.~V.~Troitsky \cite{feltro}, if $G$ is an abelian group, then the twisted conjugacy class $[e]_{\varphi}$ of the unit element $e$ is a subgroup of group $G$ for every automorphism  $\varphi \in {\rm Aut}~G$ and other twisted conjugacy classes are cosets of this subgroup. Then there arises another question if the converse statement holds.

In this paper we show, that in general case the answer to this question is negative. We show that for every nilpotent group of nilpotency class 2 the twisted conjugacy class $[e]_{\varphi}$ of the unit element is a subgroup for any automorphism $\varphi$, which is acting identicaly by modulo derived subgroup.
A.~E.~Zalessky \cite{Z} constructed the example of the infinitely generated nilpotent group of nilpotency class 2, which has only inner automorphisms. Therefore in this group the set $[e]_{\varphi}$ is a subgroup for any automorphism $\varphi \in {\rm Aut}~G,$
but $G$ is not abelian.

 We expect that the following conjecture holds.

{\bf Conjecture 1.} If $G$ is such a group where the $\varphi$-conjugacy class $[e]_{\varphi}$ of the unit element $e$ is a subgroup for every automorphism $\varphi \in {\rm Aut}~G$, then the group $G$ is nilpotent. If also $G$ is the finitely generated group, then $G$ is abelian.

As mentioned above, if $G$ is an abelian group, then for every automorphism $\varphi$ we have
$[e]_{\varphi}\leq G$ and $R(\varphi)=|G:[e]_{\varphi}|$.
In this paper we show, that even if $[e]_{\varphi}$ is a subgroup of the group $G$, then the equality $R(\varphi)=|G:[e]_{\varphi}|$  must not be fulfilled.

In the paper of E.~G.~Kukina and V.~A.~Roman'kov \cite{romankukina} there is an another formula for computation of the number $R(\varphi)$ in the torsion-free nilpotent groups. In this paper we show that this formula does not hold for the finite nilpotent groups.

Now we describe the structure of the article. In \S~1 we state general properties of the twisted conjugacy classes of the unit element. Particularly, we show that if the set $[e]_{\varphi}$ is a subgroup for any endomorphism $\varphi \in \mathrm{End}\,G$, then it is a normal subgroup of the group $G$ (Proposition 1). Example 1 shows that if two classess $[e]_{\varphi}$ and $[e]_{\psi}$ are subgroups for some automorphisms $\varphi, \psi$, then the class $[e]_{\varphi \psi}$ is not necessary a subgroup.

 In \S~2 we study the twisted conjugacy classes of the unit element for inner automorphisms. In the Proposition 4 we state the connection between the twisted conjugacy classes and  conjugacy classes. Then (Theorem 1)  we prove, that if $G$ is such a group where the twisted conjugacy class of the unit element is a subgroup for every inner automorphism, then $G$ belongs to the
Kurosh-Chernikov class $\overline{Z}.$
 In the Theorem 2 we construct a strictly descending series of the normal subgroups for such groups, where the twisted conjugacy class of the unit element is a subgroup for any inner automorphism. Thus we state that if the group has no center and satisfies the descending chain condition, then there exists such an inner automorphism of this group, that the twisted conjugacy class of the unit element is not a subgroup.
In Theorem 3 we state that if the twisted conjugacy class of the unit element is a subgroup of the group $G$ for every inner automorphism, and also $G$ satisfies the descending and ascending chain conditions, then this group is nilpotent. Therefore the first part of our conjecture holds for such groups and particularly for finite groups.

In \S~3 we prove that for every non-abelian simple group there exists such an inner automorphism $\varphi$ of this group, that the class $[e]_{\varphi}$ is not a subgroup. Then we study permutation group $S_n$ and prove, that for every inner automorphism $\varphi$, which is induced by conjugation of even permutation, the class $[e]_{\varphi}$ is not a subgroup.

 In \S~4 we study the twisted conjugacy class of the unit element in some solvable and nilpotent groups. The twisted conjugacy class $[e]_{\varphi}$ of the unit element is stated to be a subgroup for any central automorphism $\varphi$. As a consequence, we show that for any nilpotent group of nilpotency class 2 the set $[e]_{\varphi}$ is a subgroup for every IA-automorphism $\varphi$, but this statement does not hold for nilpotent groups of nilpotency class 3.

There are some open questions throughout the paper.

The authors are gratefull to O.~V.~Bruhanov, V.~A.~Roman'kov, Yu.~V.~Sosnovsky and all the participants of the seminar ``Evariste Galois'' for their helpfull comments and valuable suggestions.

\vspace{0.8cm}

\begin{center}

{\bf \S~1. Definitions and basic properties}

\end{center}

\vspace{0.5cm}

Let $G$ be a group, $\mathrm{End}\,G$ be the set of its endomorphisms, $\mathrm{Aut}\,G$ be the set of its automorphisms. If $h\in G$, then the symbol $\widehat{h}$ denotes the inner automorphism $\widehat{h}: x \mapsto h^{-1}xh$.
 We use the symbol $\mathrm{Inn}\,G$ to denote the group of all inner automorphisms. The commutator of two elements $x, y \in G$ is
$[x, y] = x^{-1} y^{-1} x y$. The following commutator identities hold in every group.
$$
[x y, z] = [x, z]^y~ [y, z] = [x, z]~[[x, z], y] ~[y, z], \eqno{(1)}
$$
$$
[x, y z] = [x, z] ~[x, y]^z = [x, z]~ [x, y] ~[[x, y], z], \eqno{(2)}
$$
$$
[x, y]^{-1} = [y, x]. \eqno{(3)}
$$
We use a symbol $g^G$ to denote a set of elements in $G$, which are conjugated with the element $g$.

If $\varphi \in \mathrm{End}\,G$, then we say, that two elements $x$ and $y$ from $G$ are
$\varphi$-{\it conjugated}, if the equality
$$
y=z^{-1}xz^{\varphi}
$$
holds for some element $z\in G$.
It is easy to see that the relation of $\varphi$-conjugation is equivalence relation on the group $G$, then the equivalence class of the element $x\in G$ is denoted $[x]_{\varphi}$. Thus
$$
[x]_{\varphi} = \{ z^{-1} x z^{\varphi}~|~z \in G \}.
$$
In particular, if $\varphi = \widehat{h}$ is an inner automorphism, we write $[x]_{h}$ instead of $[x]_{\widehat{h}}$ for the sake of simplicity. Then
$$
[x]_{h} = \{ z^{-1} x z^{h} = z^{-1} x h^{-1} z h~|~z \in G \}.
$$

It also follows from the definition that if $\varepsilon$ is the trivial endomorphism of group $G$, i.~e.
$\varepsilon (G)=e$, then
$[e]_{\varepsilon}=G$. On the other hand, for the identical automorphism $\mathrm{id}$
of the group $G$ the class $[e]_{\mathrm{id}}$ contains only the unit element $e$.
The next proposition is more interesting.

{\bf Proposition 1.} {\it Let $\varphi \in \mathrm{End}\, G$.
If  $\varphi$-conjugacy class $[e]_{\varphi}$ of the unit element $e$ is a subgroup of the group $G$, then this subgroup is normal in $G$.}

{\bf Proof.} By definition every element of the class $[e]_{\varphi}$ has the form $z^{-1}z^{\varphi}$ for some $z$ in $G$. Rewrite it in the form
$[z,\varphi]$. Then we use the commutator identity (1):
$$
[xy,\varphi]=[x,\varphi]^y[y,\varphi].
$$
We have
$$
[x,\varphi]^y=[xy,\varphi][y,\varphi]^{-1}.
$$
Since $[e]_{\varphi}$ is a subgroup of the group $G$, then $[x,\varphi]^y\in [e]_{\varphi}$ for every $y$ in $G$.
Proposition is proved.

\vspace{0.5cm}

{\bf Proposition 2.} {\it Let $\varphi\in \mathrm{End}\, G$, $\theta \in \mathrm{Aut}\, G$.
Then}

1) $([e]_{\varphi})^{\theta}=[e]_{\varphi^{\theta}}$;

2) {\it The intersection $\bigcap\limits_{\psi \in \mathrm{Aut}\, G} [e]_{\psi} $
is a characteristic subgroup of the group} $G$.

{\bf Proof.} 1) It follows directly from the following equality
$$
[x,\varphi]^\theta=[x^\theta,\varphi^\theta].
$$
2) It follows from 1). Proposition is proved.
\vspace{0.5cm}

{\bf Proposition 3.} {\it Let $\varphi,\, \psi \in \mathrm{Aut}\, G$ and
$[e]_{\varphi},\,\,\, [e]_{\psi} \leq G$.
Then $[e]_{\varphi \psi}\subseteq [e]_{\varphi}[e]_{\psi}$.}

{\bf Proof} follows from the commutator identity
$$
[x,\varphi \psi]=[x, \psi] [x,\varphi ]^\psi=[x, \psi] [x,\varphi ][[x,\varphi ],\psi]
$$
and normality of the subgroups $[e]_{\varphi}$, $[e]_{\psi}.$

\vspace{0.5cm}

The following example shows that in general case the inclusion
$[e]_{\varphi \psi}\subseteq [e]_{\varphi}[e]_{\psi}$ is strict.
 Moreover, this example shows that, if
$[e]_{\varphi}$ and $[e]_{\psi}$ are subgroups of the group $G$ for some automorphisms $\varphi$ and $\psi$, then the class $[e]_{\varphi \psi}$ must not be a subgroup.

{\bf Example 1.}
Consider free two generated nilpotent group $N = N_{2,2}$ of nilpotency class 2, with free generators $x$, $y$ and two of its automorphisms
$$
\varphi: x \mapsto x^{-1},\,\,\, y\mapsto y
$$
and
$$
\psi: x \mapsto x,\,\,\, y\mapsto y^{-1}.
$$
We will show that $[e]_{\varphi}$, $[e]_{\psi} $ are subgroups of $N$, but $[e]_{\varphi \psi}$
is not a subgroup.

Every element $z$ of $N$ can be uniquely written in the form
$$
z=x^a y^b [x,y]^c,\,\,\, a,b,c \in \mathbb{Z}.
$$
Therefore
$$
z^{-1}z^\varphi=x^{-2a} [x,y]^{-2ab-2c}.
$$
Since $a$, $b$, $c$ are arbitrary integers then
$$
[e]_{\varphi}=\left\{ \, x^{2a} [x,y]^{2c} \, |   \, a,c \in \mathbb{Z}\,  \right\}.
$$
Similarly
$$
[e]_{\psi}=\left\{ \, y^{2b} [x,y]^{2c} \, |   \, b,c \in \mathbb{Z}\,  \right\}.
$$
 $[e]_{\varphi}$, $[e]_{\psi} $ are obvious to be the subgroups of the group $N$.

But for their product $\varphi\psi $ we have
$$
z^{-1}z^{\varphi\psi}=x^{-2a} y^{-2b}[x,y]^{-2ab}.
$$
Since $a$, $b$ are arbitrary integers then
$$
[e]_{\varphi \psi}=\left\{ \, x^{2a} y^{2b}[x,y]^{-2ab} \, |   \, a,b \in \mathbb{Z}\,  \right\}.
$$
The elements $x^2$, $y^2$ belong to $[e]_{\varphi \psi}$, but their product $x^2y^2$ does not belong to $[e]_{\varphi \psi}$, therefore  $[e]_{\varphi \psi}$
is not a subgroup.

Now we show that the Fel'shtyn-Troitsky theorem does not hold for permutation group $S_3$, which is a semidirect product of cyclic group of  order 3 by a cyclic group of order 2.

{\bf Example 2.}
The permutation group on three symbols is
$$
S_3=\{e, (12), (13), (23), (123), (132) \}.
$$
We consider an inner automorphism induced by the element $h=(123)$. In this case
$$
[e]_h = [S_3, (123)] = \{ e, (132) \},
$$
i.~e. $[e]_h$ is not a subgroup.

But if we consider the inner automorphism induced by the element $h=(12),$ then we have $$
[e]_h = [S_3, (12)] = \{ e, (123), (132) \},
$$
i.~e. $[S_3,(12)]=A_3$, therefore $|S_3 : [e]_h| = 2.$ On the other hand, it is easy to see that the number $R(\widehat{h})$ of $\widehat{h}$-conjugacy classes in $S_3$ is equal to 3, as
$$
S_3 = [e]_h \sqcup [(12)]_h \sqcup [(13)]_h,
$$
where
$$
[(12)]_h = \{ (12) \},~~~ [(13)]_h = \{ (13), (23) \}.
$$

\vspace{0.8cm}

\begin{center}

{\bf \S~2. Inner automorphisms}

\end{center}

\vspace{0.5cm}
Now we show that the investigation of the twisted conjugacy classes for inner automorphisms can be reduced to the investigation of conjugacy classes.

\vspace{0.5cm}

{\bf Proposition 4.} {\it Let $g, h \in G$. Then the $\widehat{h}$-conjugacy class of the element $g$ is equal to the conjugacy class of the element $g h^{-1}$ multiplied by the element $h$, i.~e.}
$$
[g]_h = (g h^{-1})^G \cdot h.
$$

The proof follows directly from the definition.

\vspace{0.5cm}

 The number of elements conjugated to the element $g$ in $G$ is well known ( \cite[Theorem 2.5.6]{OTG}) to be equal to $|g^G| = |G : N_G (g)|$. This equality implies the following statement.

\vspace{0.5cm}

{\bf Corollary 1.} {\it The equality
$$
|[g]_h| = |G : N_G (g h^{-1})|
$$
holds for arbitrary elements $g, h \in G$.}

\vspace{0.5cm}

Particularly, if we consider $\widehat{h}$-conjugacy class of the unit element, then
$$
[e]_h = (h^{-1})^G h = [G, h].
$$
And, consequently, the set $[e]_h$ is a subgroup if and only if the set of the commutators $[G, h]$
is a subgroup. Particularly, for arbitrary elements $x, y$ in $G$ there exists an element $z \in G$ such that the equality
$$
[x, h] [y, h] = [z, h]
$$
holds, and for an arbitrary element $x \in G $ there exists an element $x' \in G$ such that the equality $$
[x, h]^{-1} = [x', h]
$$
holds.
In particular, it implies that the examples of the groups, for which the set $[e]_h$ is a subgroup for every inner automorphism $\widehat{h}$, can be found in the set of such groups, where every element of the derived subgroup of this group is a commutator. The examples of such groups are symmetrical group and alternating group. We study these groups more detailed in the next paragraph.

The next lemma is obvious.

\vspace{0.5cm}

{\bf Lemma 1.}
{\it Class $[e]_{h},$ does not contain an element $h$  if $e\neq h \in G.$}

\vspace{0.5cm}

We remind ( \cite[\S~22]{OTG}) the group belongs to Kurosh-Chernikov class $\overline{Z}$ if every series in this group can be impacted to central series. We formulate the
\vspace{0.5cm}

{\bf Theorem 1.}
{\it Let $G$ be such a group where for every element $h\in G$ the twisted conjugacy class $[e]_{h}$ of the unit element is a subgroup of $G$. Then $G$ belongs to the class $\overline{Z}$.}

{\bf Proof.} In the book of D.~J.~S.~Robinson \cite{R} there is the following criterion of belonging to the class $\overline{Z}$  for the group $N$: $N$ belongs to $\overline{Z}$ if and only if the condition
$$
h \not \in \langle [g, h] ~|~ g \in N \rangle
$$
holds for every element $h \in N$.
If $G$ satisfies the conditions of the theorem, then, as it was noted above
$$
[e]_h = [G, h] = \langle [g, h] ~|~ g \in G \rangle
$$
and by the Lemma 1 the element $h$ does not belong to $[e]_h$. And therefore $G$ belongs to the class $\overline{Z}$.

\vspace{0.5cm}

{\bf Proposition 5.}
{\it Let $G$ be such a group, where for any element $h\in G$ the twisted conjugacy class $[e]_{h}$ is a subgroup of $G$. Then for every normal subgroup $N\trianglelefteq G$ the twisted conjugacy class
$[\overline{e}]_{\overline{h}}$ is a subgroup in $G/N$.
Here $\overline{e}$ is the unit element of the factor-group $G/N$,
$\overline{h}=hN\in G/N$.}

{\bf Proof}
follows from the equality:
$$
[\overline{x},\overline{h}]~[\overline{y}, \overline{h}]=
[x,h ]~[y, h]~N =[z,h]~N=[\overline{z},\overline{h}],
$$
$$
[\overline{x},\overline{h}]^{-1}=
[x,h ]^{-1}~N =[x',h]~N=[\overline{x'},\overline{h}].
$$

\vspace{0.5cm}

{\bf Proposition 6.}
{\it Let $G$ be such a group, where for every $h\in G$ the twisted conjugacy class
$[e]_{h}$ is a subgroup of the group $G$. Then for every maximal normal proper subgroup $N$ of group $G$ the factor $G/N$ is a cyclic group of prime order.}

{\bf Proof.}
 Since $N$ is a maximal normal subgroup, then $G/N$ is a simple group. By the Proposition 5 the class $[\overline{e}]_{\overline{h}}$ is a subgroup of $G/N$ for every element $\overline{h}=hN\in G/N$.
By the Corollary 1 the factor $G/N$ is abelian. It is obvious that the simple abelian group is a cyclic group of prime order.
Proposition is proved.

\vspace{0.5cm}

Now we are ready to prove the main result about the structure of such a group, where the twisted conjugacy class of the unit element is a subgroup for every inner automorphism.

{\bf Theorem 2.} {\it Let $G$ be such a group, where for every element $h \in G$ the $\widehat{h}$-conjugacy class
$[e]_{h}$  is a subgroup. Then there is a strictly descending series of normal subgroups in $G$}:
$$
G \rhd G_1 \rhd G_2 \rhd \ldots
$$

{\bf Proof. } Choose some element $h_0\in G$ and define the subgroup $G_1=[e]_{h_0}$. By the Lemma 1 we have the inequality $G_1\neq G$ as $h_0$ does not belong to $G_1$.
Then choose some element $h_1\in G_1$ and define the subgroup $G_2=[e]_{h_1}$.
As $G_1$ is a normal subgroup of  $G$, then $G_2$ is contained in $G_1$. By the Lemma 1 we have $G_2\neq G_1$ as $h_1$ belongs to $G_2$, and so on. We receive the required series of the normal subgroups in group $G$:
$$
G\rhd G_1=[e]_{h_0} \rhd G_2=[e]_{h_1} \rhd ... \rhd G_{k+1}=[e]_{h_k} \rhd ...,
$$
where $h_k \in G_k=[e]_{h_{k-1}}$, $k=1, 2, ...$ The theorem is proved.

\vspace{0.5cm}

{\bf Question 1.} What can we say about the factors of this series?

The following statement follows from this theorem.

\vspace{0.5cm}

{\bf Corollary 2.}
{\it If the group $G$ has no center and it satisfies the descending chain condition for normal subgroup, then there exists an element $h\in G,$ such  that the twisted conjugacy class $[e]_{h} $ is not a subgroup of} $G$.

{\bf Proof.} Suppose the contrary. Let $G$ be such a group, that for every element $h\in G$ the twisted conjugacy class
$[e]_{h}$ is a subgroup of $G$. By the Theorem 2 we have a normal series:
$$
G\rhd G_1=[e]_{h_0} \rhd G_2=[e]_{h_1} \rhd ... \rhd G_{k+1}=[e]_{h_k} \rhd ...,
$$
where $e\neq h_k \in G_k=[e]_{h_{k-1}}$, $k=1, 2, ...$
From the descending chain condition it follows that
$ G_{k+1}=\{e\}$ for some $k$. But it means that for every element $x\in G$
$$
x^{-1}h_k^{-1}xh_k=e,
$$
i.~e. $h_k$ belongs to the center of group $G$. But as $h_k \neq e$, then it contradicts the condition that the group $G$ has no center.
Corollary is proved.

The following proposition follows directly from the proof of this corollary.

\textbf{Proposition 7.} {\it Let $[e]_h$ be a subgroup for some element $h \in G$. Then the center of the factor group $G/[e]_h$ is not trivial.}

\textbf{Proof.} The element $h [e]_h \neq [e]_h$ is obvious to belong to the center of its group.

\vspace{0.5cm}

Later we use the symbol min-n (max-n) to denote the descending (respectively ascending) chain condition of normal subgroups. The following statement holds.

{\bf Corollary 3.}
{\it Let $G$ satisfies min-n condition and for every element $h \in G$ class $[e]_h$ is a subgroup of $G$.
Then the center of the group $G$ is not trivial.}

{\bf Proof.} By the Theorem 2 there exists a normal series
$$
G\rhd G_1=[e]_{h_0} \rhd G_2=[e]_{h_1} \rhd ... \rhd G_{k+1}=[e]_{h_k} \rhd ...
$$
According to the min-n condition it follow that this series is finite, i.~e. $[e]_{h_n} = \{ e \}$ for some nontrivial element $h_n$, but it means that the element $h_n$ belongs to the center of group $G$. Theorem is proved.

\vspace{0.5cm}

{\bf Theorem 3.}
{\it Let a group $G$ satisfies the min-n and max-n conditions, and for every element $h \in G$ the class $[e]_h$ is a subgroup of the group $G$.
Then $G$ is nilpotent group.}

{\bf Proof.} By the Corollary 3 the center $Z(G)$ of the group $G$ is not trivial. By the Proposition 5 the factor group $\overline{G} = G / Z(G)$ satisfies the condition that for any element $\overline{h} \in \overline{G}$ class $[\overline{e}]_{\overline{h}}$ is a subgroup of $\overline{G}$. Taking into account the fact, that the min-n condition holds for the factor groups, we conclude that the center of group $\overline{G}$ is not trivial.
Therefore, there is an ascending series of hypercenters in  $G$
$$
\{ e \} < \zeta_1(G) = Z(G) < \zeta_2(G) < \ldots
$$
Since $G$ satisfies the max-n condition, we conclude that $\zeta_m(G) = G$ for some positive integer $m$.
And consequently, $G$ is nilpotent group. The theorem is proved.

\vspace{0.5cm}

Every finite group satisfies the min-n and max-n conditions, ant therefore from the Theorem 3 we have

{\bf Corollary 4.}
{\it Let  the twisted conjugacy classes $[e]_h$ is a subgroup of the finite group $G$ for every $h \in G$.
Then $G$ is nilpotent group.}

\vspace{0.5cm}

\textbf{Question 2.}  Can we state that in this case $G$ is nilpotent group of nilpotency class 2?

\vspace{0.5cm}

Now we determine the relation between twisted conjugacy classes of the unit element and some verbal subgroups.

{\bf Proposition 8.}
{\it Let the twisted conjugacy class $[e]_{h}$ be a subgroup of group $G$ for every element $h\in G$.
If for any nontrivial class $[e]_{h} \neq \{e\}$ the factor group $G/[e]_{h}$ is abelian, then the derived subgroup $G'$ belongs to the center $Z(G)$ of the group $G$.}

{\bf Proof.}
Since $G/[e]_{h}$ is abelian, then  $G'\subseteq [e]_{h}$.
On the other hand $[e]_{h}\subseteq G'$, and therefore $G'= [e]_{h}$.
Let $h_1\in G'$ and $h_1\not\in Z(G)$.
Then $[e]_{h_1} < [e]_{h}$ (the inclusion is strict since $h_1\not\in [e]_{h_1}$).
It contradicts to the fact that $[e]_{h_1} = G' = [e]_{h}$.
The proposition is proved.

Let now $w(x) = w(x_1,\dots,x_n)$ be an outer commutator word in the free group
$F_n = \langle x_1, x_2, \ldots, x_n \rangle$, i.~e. there are variables $x_1,x_2, \dots, x_n$ and $n-1$ pairs of commutator brackets between it. If $G$ is a group, then we use the symbol $w(g) = w(g_1, g_2, \ldots, g_n)$ to denote the value of $w$ in $G$. In these notations we have a

\textbf{Lemma 2.} {\it Let  the twisted conjugacy class $[e]_{h}$ is a subgroup of the group $G$ for every element $h\in G$, and $g = (g_1, g_2, \ldots, g_n)\in G^n.$}
Then

1) $w(g) \in [e]_{g_1}\cap \dots \cap [e]_{g_n};$

2) $[e]_{w(g_1, \dots, g_n)}\leq [e]_{g_1}\cap \dots \cap [e]_{g_n}$.

\textbf{Proof.} 1) We use the induction by the length of the word $w(x)$. If $n=2$ then $w(g)=[g_1,g_2]\in [e]_{g_2}$.
On the other hand $w(g)=[g_2,g_1]^{-1}\in [e]_{g_1}$.  Assume that the statement holds for all the words, which are shorter than $n$. Let $w=[u,v]$ be a commutator word of the length $n$, where $u,v$ are commutator word of the less length. For the definiteness we assume that $u=u(x_1,\dots,x_k),v=v(x_{k+1},\dots,x_n)$.
Let $i\leq k$. By the induction hypothesis $u(g) \in [e]_{g_i}$. Then $u(g)=[a_i, g_i]$
for any $a_i$, $i=1, \ldots,k$. Consequently $w(g)=[[a_i, g_i],v(g)]=[a_i, g_i]^{-1}[a_i, g_i]^v
\in [e]_{g_i}$. Since $i$ is an arbitrary number we have $w(g) \in [e]_{g_1} \cap \dots \cap [e]_{g_k}$.
The case $i\geq k$ is similar.

2) Since $w\in [e]_{g_i}$ and $[e]_{g_i}$ is normal subgroup of $G$, then $[e]_{w(g)}$ is a subgroup of $[e]_{g_i}$.  Lemma is proved.

\vspace{0.5cm}

In the theorem 2 we construct a normal series of subgroups in $G$. The following proposition shows what elements can be chosen as $\{h_i \}.$

\vspace{0.5cm}

\textbf{Proposition 9.} {\it For every element $g_1, \dots, g_n$ in $G$ there exists a normal series
$$
[e]_{g_1}>[e]_{[g_1, g_2]}>[e]_{[g_1, g_2, g_3]}> \dots > [e]_{[g_1, \dots, g_n]}.
$$ Moreover, if all element
$${g_1},{[g_1, g_2]},{[g_1, g_2, g_3]}, \dots ,
{[g_1,\dots, g_n]}
$$
are not central, then all inclusions are strict.}

\textbf{Proof.} If $a=[g_1, \dots, g_i], b=[g_1, \dots, g_{i+1}]$, then $b=[a,g_{i+1}]$ and by the Lemma 2 $$[e]_{[g_1,\dots,g_{i+1}]}=[e]_{[a,g_{i+1}]}=[e]_b\leq [e]_{a}=[e]_{[g_1,\dots,g_i]}.$$
If $[g_1, \dots, g_i+1]$ does not belong to $Z(G)$, then $[g_1, \dots, g_{i+1}]\in[e]_{[g_1,\dots,g_{i}]}\setminus[e]_{[g_1,\dots,g_{i+1}]}. $
Proposition is proved.

\vspace{0.5cm}

As already noted, the intersection $\bigcap\limits_{\varphi \in \mathrm{Aut}\, G} [e]_{\varphi}$ is a characteristic subgroup. Moreover, by the Lemma 1 we have
$$\bigcap\limits_{\varphi \in \mathrm{Aut}\, G}
[e]_{\varphi}<\bigcap\limits_{g \in G\backslash Z(G)} [e]_{g}\subset Z(G).$$

%\textbf{Предложение .} Пусть $G$ группа, a $\varphi : G\rightarrow G$ -- центральный автоморфизм группы $G$.
%Тогда множество $[e]_{\varphi}$ является подгруппой в $G$.

%\textbf{Доказательство.} Так как $\varphi$ -- центральный автоморфизм группы $G$, то он действует по правилу
%$$\varphi: x \mapsto \overline{\varphi}(x)x,$$
%где $\overline{\varphi}:G\rightarrow Z(G)$ -- гомоморфизм группы $G$ в её центр $Z(G)$. Отсюда следует,
%что $$[e]_{\varphi}=\{x\varphi(x^{-1})~|~ x\in G\}=\{x\overline{\varphi}(x^{-1})x^{-1}~| ~ x\in G\}=\{\overline{\varphi}(x)~|~ x\in G\}$$
%Очевидно, что множество такого вида является подгруппой в $G$. Более того, так как
%$\overline{\varphi}$ -- гомоморфизм в центр группы $G$, то $[e]_{\varphi}$ является подгруппой и в $Z(G)$.
%\begin{flushright}
%$\square$
%\end{flushright}

\vspace{0.8cm}

\begin{center}

{\bf \S~3. Simple groups and permutation groups}

\end{center}

\vspace{0.5cm}

For simple groups we have

{\bf Theorem 4.}
{\it Let $G$ be a simple non-abelian group and $h$ be a nontrivial element of $G$. Then $[e]_{h}$ is not a subgroup of $G$.}

{\bf Proof.}
Suppose the contrary.
Let $[e]_{h}$ be a subgroup of the group $G$.
By the Proposition 1 this subgroup is normal.
As $G$ is a simple group, then $[e]_{h}=G$.
But it means that for every element $g\in G$ there exists an element $x\in G$, for which the equality
$$
[x, h] = g
$$
holds.
If $g=h$ we have $[x, h] = h$, i.~e. $h = e$.
It contradicts to the fact that $h$ is nontrivial element of group $G$.
Theorem is proved.

\vspace{0.5cm}

{\bf Remark.}
If class $[e]_{h}$ is a subgroup, and by the Proposition 1 a normal subgroup of the group $G$, then by the Proposition 2 the equality $([e]_{h})^g=[e]_{h^g}$ holds for every $g\in G$. Therefore, going through the element $h$ in $G$, it is sufficient to go through the transversal set of conjugacy classes.

\vspace{0.5cm}

For the alternating group we have a

{\bf Proposition 10.} {\it In the group $A_n$, $n \geq 4$ there exists an element $h$, such that the class $[e]_h$ is not a subgroup.}

{\bf Proof.}
If $n \geq 5$, then this statement follows from the Theorem 4.

Consider the group
$$
  A_4= \{e, (123), (132), (134), (143), (124), (142), (234), (243), (12)(34), (13)(24), (14)(23)\}.
$$
By the remark, which is stated above, we have to consider two cases
$$
h = (123),\,\,\,  (12)(34).
$$

 Let $h=(123)$ and consider the class $[e]_{h}$. It is obvious that it coincides with the set of commutators $[A_4, h]$. By the direct calculations we have
$$
[A_4, (123)]= \{ e, (12)(34), (13)(24), (14)(23) \},
$$
i.~e.
$[e]_{(123)}$
is a Klein's four-group \cite[\S~3]{OTG}.
In this case the element $(123)$ does not belong to this group.

 Let now $h=(12)(34)$. Then
$$
[A_4, (12)(34)]= \{ e, (13)(24), (14)(23) \}.
$$
It is obvious that these elements do not form a subgroup, and therefore the set $[e]_{(12)(34)}$ is not a subgroup.
Proposition is proved.

\vspace{0.5cm}

Consider now the permutation group $S_n.$ The following statement holds.

{\bf Proposition 11.} {\it Let the element $h \in A_n$ be in the permutation group $S_n$, $n \geq 5$. Then the class $[e]_h$ is not a subgroup.}

{\bf Proof.} Since $[e]_h = [S_n, h]$, all the elements of $[e]_h$ belong to $A_n$. If $[e]_h$ is a subgroup, then by the Proposition 1 this subgroup is normal. Since the center of the group $S_n$ is trivial, we conclude that $[e]_h$ contains a non-trivial element, and therefore it coincides with $A_n$. On the other hand, by the Lemma 1 the element $h$ does not belong to the class $[e]_h$. This contradiction proves the proposition.

 Note that groups $S_3$ and $A_4$ are semidirect products of the abelian group by the cyclic group. Really, there is a shot exact sequence
$$
1 \longrightarrow A_3 \longrightarrow S_3 \longrightarrow \mathbb{Z}_2 \longrightarrow 1,
$$
   and $A_3=S_3'$ and $S_3=\mathbb{Z}_3 \leftthreetimes \mathbb{Z}_2$, i.~e $S_3$ is an extension of the group $\mathbb{Z}_3$ by the group $\mathbb{Z}_2$.

Similarly there is an exact sequence
$$
1 \longrightarrow K_4 \longrightarrow A_4 \longrightarrow \mathbb{Z}_3 \longrightarrow 1,
$$
and $K_4=A_4'=\mathbb{Z}_2 \times \mathbb{Z}_2$, i.~e $A_4$ is an extension of the group $\mathbb{Z}_2 \times \mathbb{Z}_2$ by the group $\mathbb{Z}_3$.

\vspace{0.5cm}

{\bf Questions.}
4. Let $G$ be an extension of the  abelian group by the cyclic (abelian) group. Can we state that in this case there is such an element $h \in G$, for which the class $[e]_h$ is not a subgroup?

5. The extensions of the abelian groups are described in terms of homology groups. Is it possible to formulate in this terms a condition that every twisted conjugacy class of the unit element is a subgroup of the group?

\vspace{0.8cm}

\begin{center}

{\bf \S~4. Solvable groups}

\end{center}

\vspace{0.5cm}

As shown in the Example 1, there is  a nilpotent group of nilpotency class 2 and an automorphism of this group, such that the twisted conjugacy class of the unit element is not a subgroup. However for central automorphisms we have

\textbf{Proposition 12.} {\it Let $G$ be a group and $\varphi : G\rightarrow G$ a central of $G$. Then the set $[e]_{\varphi}$ is a subgroup of $G$.}

\textbf{Proof.} Since $\varphi$ is a central automorphism of the group $G$, then it acts by the following way
$$
\varphi: x \mapsto \overline{\varphi}(x)x,
$$
 where $\overline{\varphi}:G\rightarrow Z(G)$ is a homomorphism of the group $G$ into its center $Z(G)$. Therefore
$$
[e]_{\varphi}=\{x^{-1} \varphi(x)~|~ x\in G\}=\{x^{-1} \overline{\varphi}(x)x~| ~ x\in G\}
=\{\overline{\varphi}(x)~|~ x\in G\}.
$$
This set is obvious to be a subgroup of the group $G$. Moreover, since $\overline{\varphi}$ is a homomorphism into the center of the group $G$, then $[e]_{\varphi}$ is a subgroup of $Z(G)$.

\medskip

{\bf Conjecture 2.}
If the group $G$ has no center, then  there exists an inner automorphism $\varphi$, such that $[e]_{\varphi}$ is not a subgroup.

\vspace{0.5cm}

Thus we need to consider the groups with nontrivial center, in particular, nilpotent groups.

As we know, the group $S_n$ is perfect for $n\geq 3$, and therefore it has no center. The group $A_n$ is simple for $n\geq 5$, and therefore it has no center.
And as it shown above, the group $A_4$ also has no center.

\vspace{0.5cm}

The automorphism of the free nilpotent group $N_{n,r}$ is called an IA-{\it automorphism}, if it acts identically by modulo derived subgroup of $N_{n,r}$. There is a shot exact sequence
$$
1 \longrightarrow {\rm IA}(N_{n,r}) \longrightarrow {\rm Aut} N_{n,r} \longrightarrow {\rm GL}_n(\mathbb{Z}) \longrightarrow 1.
$$
If $r=2$ then every IA-automorphism is a central automorphism, and therefore by the proposition 2 we have a

{\bf Corollary 5.} {\it Let $\varphi$ be an IA-automorphism of the free nilpotent group $N_{n,2}$ of nilpotency class 2. Then the set $[e]_{\varphi}$ is a subgroup.}

\medskip

Now we show that this statement can not be generalized to the free nilpotent groups of greater nilpotency class.

{\bf Proposition 13.} {\it There is such an inner automorphism of the 2 generated free nilpotent group of nilpotency class 3, that the twisted conjugacy class of the unit element is not a subgroup.}

{\bf Proof.}
Let $G=N_{2,3}$ be the free nilpotent group of nilpotency class 3 with free generators $x,y$ and $\varphi$ be the following automorphism
$$
x \mapsto x[x,y], ~~y \mapsto y.
$$
This automorphism is obvious to be $IA$-automorphism. Moreover, $\varphi$ is an inner automorphism:
$\varphi = \widehat{y}$ conjugaction by the element $y$.

Every element of the group $G$ can be uniquely expressed in the form
$$
g=x^ay^b[y,x]^c[[y,x],y]^d[[y,x],x]^f,
$$ for some $a, b, c, d, f$. Then
$$
[g,y]=[x^ay^b[y,x]^c,y]=[x^ay^b,y]^{[y,x]^c}[[y,x]^c,y]=[x^ay^b,y][[y,x],y]^c=
[x^a,y]^{y^b}[[y,x],y]^c=
$$
$$
=[x^a,y][[x^a,y],y^b][[y,x],y]^c=
[x,y]^a[[x,y],x]^{\frac{a(a-1)}{2}}[[x,y],y]^{ab}[[y,x],y]^c=
$$
$$
=[x,y]^a[[x,y],x]^{\frac{a(a-1)}{2}}[[y,x],y]^{c-ab}=[y,x]^{-a}[[y,x],x]^{-\frac{a(a-1)}{2}}[[y,x],y]^{c-ab},
$$
i.~e. the twisted conjugacy class of the unit element is
$$
[e]_{\varphi}=\{[y,x]^{-a}[[y,x],x]^{-\frac{a(a-1)}{2}}[[y,x],y]^{c-ab}~|~a,b,c \in \mathbb{Z}\}.
$$
It is clear that the commutator $[x,y]$ belongs to the class $[e]_{\varphi}$ (for $a = b = c = 1$). But its square $[x,y]^2$ does not belong to the class $[e]_{\varphi}$. Really, if we suppose that $[x,y]^2\in [e]_{\varphi}$, then
$$[x,y]^2=[y,x]^{-a}[[y,x],x]^{-\frac{a(a-1)}{2}}[[y,x],y]^{c-ab}.$$
And consequently $a=2$ and $$[x,y]^2=[y,x]^{-2}=[y,x]^{-2}[[y,x],x]^{-1}[[y,x],y]^{c-2b},$$
i.~e. $[[y,x],x]^{-1}[[y,x],y]^{c-2b}=e$. But it is not true in the group $N_{2,3}$ for any integers $b$ and $c$.
As a consequence the class $[e]_{\varphi}$ is not a subgroup.

\medskip

%\newpage

{\bf Proposition 14.} {\it There exists such a group $G$, that the class $[e]_{\varphi}$ is not a subgroup of $G$ for some automorphism $\varphi$, but for every endomorphism $\psi$ with non-trivial kernel the class $[e]_{\psi}$ is a subgroup of the group $G$.}

{\bf Proof.}
Consider the group
$$
G =\mbox{gr}\left( x, y   \,\|\,  x^2=y^2=[x,y]^2=[[x,y],x]=[[x,y],y]=1 \right) .
$$
Directly from this representation we have

1. $G=\left\{ e, x, y, xy,    [x,y], x[x,y], y[x,y], xy[x,y] \right\} $,
i.~e. $|G|=8$.

2. The following equalities hold in $G$
$$
[x,y]=(xy)^2,\,\, [x,y]^{-1}=[y,x]=[x,y],\,\, (xy)^4=e,\,\, yx=(xy)^3.
$$

3. The center of the group $G$ coincides with its derived subgroup, and it is a cyclic subgroup of the second order with the generator $[x,y]=(xy)^2$.

We describe the twisted conjugacy classes
$[e]_{\varphi}$, $\varphi \in \mathrm{End}~G$.

By the proposition 2 it is sufficient to choose a transversal set of $\mathrm{End}~G$ with respect to $\mathrm{Aut}~G$.

Note that any automorphism of the group $G$ is determined by its action on the generators $x$, $y$. Since $x$, $y$ are the elements of the second order, then by the defining relations of the group $G$, we see that any map of the form
$$
x^{\varphi}=a,\,\,\, y^{\varphi}=b,
$$
where $a$ and $b$ are the elements of the second order, defines the endomorphism of the group $G$.

Since the set of all elements of the second order of the group $G$ is of the form
$$
\left\{ e, x, y,  [x,y], x[x,y], y[x,y]  \right\},
$$
then $\left|\mathrm{End}~G\right|=36$. The full information about the set $\mathrm{End}~G$ is contained in the following
table:

$$
\begin{array}{|c|c|c|c|c|c|c|}
  \hline
        & e & x & y & [x,y] & x[x,y] & y[x,y] \\
  \hline
     e  & E.0 & E.2 & E.4 & E.5 & E.2 & E.4 \\
  \hline
     x  & E.4 & E.3 & A.3 & E.1 & E.7 & A.4 \\
  \hline
     y  & E.2 & A.0 & E.3 & E.8 & A.1 & E.7 \\
  \hline
  [x,y] & E.5 & E.8 & E.1 & E.6 & E.8 & E.1 \\
  \hline
  x[x,y] & E.4 & E.7 & A.4 & E.1 & E.3 & A.3 \\
  \hline
  y[x,y] & E.2 & A.1 & E.7 & E.8 & A.2 & E.3 \\
  \hline
\end{array}
$$

In the top line there is an image of the element $x$ by the endomorphism, and in the left column there is an image of the element $y$ by the endomorphism. The letter $A$ means that this map determines the automorphism, and the letter $E$ means that this map determines the endomorphism with non-trivial kernel.
The number after the letters $A$, $E$ is the number of the conjugacy class with respect to conjugation by elements of the group $\mathrm{Aut}~G$.

For example, class $E.0$ contains only the trivial automorphism $\varepsilon$:
$\varepsilon (g)=e$ for any $g\in G$.
Class $A.0$ contains the identical automorphism $G$.
Class $E.1$ contains four automorphisms
$$
(x^{\varphi}=(xy)^2,\,\,\, y^{\varphi}=x),\,\,\,
(x^{\varphi}=y,\,\,\, y^{\varphi}=(xy)^2),
$$
$$
(x^{\varphi}=y(xy)^2,\,\,\, y^{\varphi}=(xy)^2),\,\,\,
(x^{\varphi}=(xy)^2,\,\,\, y^{\varphi}=x(xy)^2).
$$
And so on. We note that the class $A.2$ consists of the unique inner automorphism, which corresponds to the element $xy$ of the group $G$, and this automorphism generates a center of the group $\mathrm{Aut}~G$. The order of the group $\mathrm{Aut}~G$ is equal to eight.

 There is the full information about all conjugacy classes for automorphisms and endomorphisms in the following tables

$$
\begin{array}{|c|c|c|c|c|c|c|}
  \hline
        \mathrm{Aut}\, G            & x^{\varphi}=x & x^{\varphi}=y & x^{\varphi}=y[x,y] & x^{\varphi}=x      & x^{\varphi}=x[x,y]  \\
                                    & y^{\varphi}=y & y^{\varphi}=x & y^{\varphi}=x      & y^{\varphi}=y[x,y] & y^{\varphi}=y[x,y]  \\
  \hline
     x^{-1}x^{\varphi}              & e             & xy            & (xy)^3             & e                  & (xy)^2              \\
  \hline
     (xy)^{-1}(xy)^{\varphi}        & e             & (xy)^2        & e                  & (xy)^2             & e                   \\
  \hline
     y^{-1}y^{\varphi}              & e             & (xy)^3        & (xy)^3             & (xy)^2             & (xy)^2              \\
  \hline
    [x,y]^{-1}[x,y]^{\varphi}       & e             & e             & e                  & e                  & e                   \\
  \hline
  (x[x,y])^{-1}(x[x,y])^{\varphi}   & e             & xy            & (xy)^3             & e                  & (xy)^2              \\
  \hline
  (xy[x,y])^{-1}(xy[x,y])^{\varphi} & e             & (xy)^2        & e                  & (xy)^2             & e                   \\
  \hline
  (y[x,y])^{-1}(y[x,y])^{\varphi}   & e             & (xy)^3        & (xy)^3             & (xy)^2             & (xy)^2              \\
  \hline
                                    & \mbox{gr}(e)  & \mbox{gr}(xy) & -                  & \mbox{gr}((xy)^2)  & \mbox{gr}((xy)^2)   \\
  \hline
\end{array}
$$

{\footnotesize

$$
\begin{array}{|c|c|c|c|c|c|c|c|c|}
  \hline
        \mathrm{End}\, G            & x^{\varphi}=x      & x^{\varphi}=x & x^{\varphi}=(xy)^2 & x^{\varphi}=(xy)^2 & x^{\varphi}=y & x^{\varphi}=y      & x^{\varphi}=x(xy)^2 &  x^{\varphi}=x(xy)^2\\
                                    & y^{\varphi}=e      & y^{\varphi}=x & y^{\varphi}=(xy)^2 & y^{\varphi}=e      & y^{\varphi}=e & y^{\varphi}=(xy)^2 & y^{\varphi}=x       &  y^{\varphi}=y      \\
  \hline
     x^{-1}x^{\varphi}              & e                  & e             & x(xy)^2            & x(xy)^2            & xy            & xy                 & (xy)^2              &  x(xy)^2         \\
  \hline
     (xy)^{-1}(xy)^{\varphi}        & y                  & (xy)^3        & (xy)^3           & xy                 & x(xy)^2       & x                  & xy                  &  x            \\
  \hline
     y^{-1}y^{\varphi}              & y                  & (xy)^3        & y(xy)^2            & y                  & y             & y(xy)^2            & (xy)^3              &  e             \\
  \hline
    [x,y]^{-1}[x,y]^{\varphi}       & (xy)^2             & (xy)^2        & (xy)^2             & (xy)^2             & (xy)^2        & (xy)^2             & (xy)^2              &  (xy)^2            \\
  \hline
  (x[x,y])^{-1}(x[x,y])^{\varphi}   & (xy)^2             & (xy)^2        & x                  & x                  & (xy)^3      & (xy)^3           & e                   &  x       \\
  \hline
  (xy[x,y])^{-1}(xy[x,y])^{\varphi} & y(xy)^2            & xy            & xy                 & (xy)^3           & x             & x(xy)^2            & (xy)^3              &  x(xy)^2                 \\
  \hline
  (y[x,y])^{-1}(y[x,y])^{\varphi}   & y(xy)^2            & xy            & y                  & y(xy)^2            & y(xy)^2       & y                  & xy                  &  (xy)^2             \\
  \hline
                                    & \mbox{gr}(y,(xy)^2)& \mbox{gr}(xy)             & G                  & G                  & G             & G                  & \mbox{gr}(xy)       & \mbox{gr}(y,(xy)^2)                             \\
  \hline
\end{array}
$$

}

So, there exists such an automorphism $\varphi$ of the group $G$, that the class $[e]_{\varphi}$ is not a subgroup, but for every endomorphism $\varphi$ from the $\mathrm{End}~G \backslash \mathrm{Aut}~G$
the class $[e]_{\varphi}$ is a subgroup.

Now we show that the equality
$$
R(\varphi) = | G : [e]_{\varphi} |
$$
does not hold for the group from the last proposition. Really, consider the automorphism $\varphi$ defined on the generators
$$
x^{\varphi} = y,~~y^{\varphi} = x.
$$
Then the $\varphi$-conjugacy class of the unit element is a subgroup, and it consists of the elements
$$
[e]_{\varphi} = \{ e, xy, (xy)^2, (xy)^3 \};
$$
the $\varphi$-conjugacy class of the element $x$ consists of the elements
$$
[x]_{\varphi} = \{ x, y \};
$$
the $\varphi$-conjugacy class of the element $x[x,y]$ consists of the elements
$$
[x[x,y]]_{\varphi} = \{  x [x, y], y [x, y] \}.
$$
Thus
$$
G = [e]_{\varphi} \sqcup [x]_{\varphi} \sqcup [x[x,y]]_{\varphi},
$$
i.~e. $R(\varphi) = 3$. On the other hand, $|G| = 8,$ $[e]_{\varphi} = 4$, and therefore $| G : [e]_{\varphi} | = 2$.

\medskip

Let $G$ be finitely generated  torsion-free nilpotent group of nilpotency class $k$ and
$$
\zeta_0G=1<\zeta_1G< \dots <\zeta_{k-1}G<\zeta_kG=G
$$ be the upper central series.
Let $A_i=\zeta_{i+1}G / \zeta_iG$, $\varphi_i$ be the automorphism induced by the automorphism $\varphi$ on $A_i$.
If $R(\varphi) < \infty$, then, as it stated in the paper \cite{romankukina}, the equality
$$
R(\varphi)=\prod_{i=0}^{k-1} |A_i:[e]_{\varphi_i}|
$$
holds.

Now we show, that this formula does not hold for finite nilpotent groups. We consider the group from the last proposition as a contrary instance. Really, for this group we have
$$
\zeta_1(G) = Z(G) = \langle (x y)^2 \rangle \simeq \mathbb{Z}_2,~~~\zeta_2(G) = G.
$$
The automorphism
$$
x^{\varphi} = y,~~y^{\varphi} = x
$$
induces the automorphisms $\varphi_0$ and $\varphi_1$ of the groups $Z(G)$ and $G / Z(G)$ respectively. And in this case
$|Z(G) : [e]_{\varphi_0}| = 2$, and $|G / Z(G) : [\overline{e}]_{\varphi_1}| = 2$. As a consequence, we have
$$
|Z(G) : [e]_{\varphi_0}| \cdot |Z(G) : [e]_{\varphi_0}| = 4,
$$
but, as we noted, $R(\varphi) = 3$.

\medskip

Now we make an example of the torsion-free group, which is an extension of the free abelian group by the infinite cyclic group, where there is such an inner automorphism, that the twisted conjugacy class of the unit element is not a subgroup.

{\bf Example 3.}
Let  $\lambda$ be a variable. We consider the group $G$ generated by the matrices
$$
 d=
 \left(%
 \begin{array}{cc}
 \lambda & 0 \\
  0 & 1 \\
 \end{array}%
 \right),\,\,\,
 t_{12}(1)=
 \left(%
 \begin{array}{cc}
  1 & 1 \\
  0 & 1 \\
 \end{array}%
 \right).
$$
The elements of the matrices belong to the ring $\mathbb{Z}[\lambda^{\pm 1}]$. As it was noted in \cite[\S~6]{OTG}, the group $G$ is isomorphic to the wreath product $\mathbb{Z} \wr \mathbb{Z}.$

Note that the formulas
$$
d^{-1} t_{12}(\mu) d = t_{12}(\lambda^{-1}\mu),\,\,\,
d t_{12}(\mu) d^{-1}=t_{12}(\lambda \mu),\,\,\,
[d^l,t_{12}(\mu)] = t_{12}(\mu(1-\lambda^{-l})),
$$
hold for $l\in \mathbb{Z}$ an for and arbitrary element $\mu$ of the ring $\mathbb{Z}[\lambda^{\pm 1}]$.

Now we show that there is such an element $h = d^m t_{12}(\mu)\in G$, that the class $[e]_{h}$ is not a subgroup of the group $G$.

Firstly, we note that the derived subgroup $G'$ of the group $G$ is an abelian group, and every element of the group $G$ can be presented in the form $h = d^m t_{12}(\mu)$.
Then by the commutator identities (1)-(3) from \S~1 we have
$$
[d^l t_{12}(\nu), d^m t_{12}(\mu)] = [d^l, d^m t_{12}(\mu)]^{t_{12}(\nu)}[t_{12}(\nu), d^m t_{12}(\mu)]=
$$
$$
 =[d^l, d^m t_{12}(\mu)][t_{12}(\nu), d^m t_{12}(\mu)]=
 [d^m t_{12}(\mu), d^l]^{-1}[d^m t_{12}(\mu), t_{12}(\nu)]^{-1}=
$$
$$
=\left([d^m, d^l]^{t_{12}(\mu)}[t_{12}(\mu), d^l]\right)^{-1}
 \left([d^m, t_{12}(\nu)]^{t_{12}(\mu)}[t_{12}(\mu),t_{12}(\nu)]\right)^{-1}=
$$
$$
 =[d^l, t_{12}(\mu)] [d^m, t_{12}(\nu)]^{-1}=
 t_{12}(\mu(1-\lambda^{-l})-\nu (1-\lambda^{-m})).
$$
And consequently,
$$
[d^l t_{12}(\nu), d^m t_{12}(\mu)][d^{l_1}t_{12}(\nu_1), d^m t_{12}(\mu)]=
 t_{12}\left(\mu(2-\lambda^{-l}-\lambda^{-l_1})-(\nu +\nu_1) (1-\lambda^{-m})\right).
$$
If
$$
[d^l t_{12}(\nu), d^m t_{12}(\mu)][d^{l_1}t_{12}(\nu_1), d^m t_{12}(\mu)]=
 [d^{l_2}t_{12}(\nu_2), d^m t_{12}(\mu)],
$$
i.~e. $[G, d^mt_{12}(\mu)]$ is a subgroup, then
$$
\mu(2-\lambda^{-l}-\lambda^{-l_1})-(\nu +\nu_1) (1-\lambda^{-m})=\mu(1-\lambda^{-l_2})-\nu_2 (1-\lambda^{-m})
$$
or
$$
\mu(1+\lambda^{-l_2}-\lambda^{-l}-\lambda^{-l_1})= (1-\lambda^{-m})(\nu +\nu_1-\nu_2).
$$
Let, for example, $m=0$ and $\mu \neq 0$. Then
$$
1+\lambda^{-l_2}=\lambda^{-l}+\lambda^{-l_1}.
$$
This equality holds if and only if
$$
 (l=0\,\,\, \mbox{and} \,\,\, l_1 = l_2)\,\,\, \mbox{or} \,\,\, (l_1=0\,\,\, \mbox{and} \,\,\,l=l_2).
$$

Thus the set $[G,d^mt_{12}(\mu)$ is not a subgroup for any $m$ and $\mu$.

%%%%%%%%%%%%%%%%%%%%%%%%%%%%%%%%%%%%%%%%%%%%%%%%%%%%%%%%%%%%%%%%%%%%%%%%%%%%%%%%%%%%%%%
%%%%%%%%%%%%%%%%%%%%%%%%%%%%%%%%%%%%%%%%%%%%%%%%%%%%%%%%%%%%%%%%%%%%%%%%%%%%%%%%%%%%%%%
%%%%%%%%%%%%%%%%%%%%%%%%%%%%%%%%%%%%%%%%%%%%%%%%%%%%%%%%%%%%%%%%%%%%%%%%%%%%%%%%%%%%%%%
%%%%%%%%%%%%%%%%%%%%%           Вопросы          %%%%%%%%%%%%%%%%%%%%%%%%%%%%%%%%%%%%%%
%%%%%%%%%%%%%%%%%%%%%%%%%%%%%%%%%%%%%%%%%%%%%%%%%%%%%%%%%%%%%%%%%%%%%%%%%%%%%%%%%%%%%%%
%%%%%%%%%%%%%%%%%%%%%%%%%%%%%%%%%%%%%%%%%%%%%%%%%%%%%%%%%%%%%%%%%%%%%%%%%%%%%%%%%%%%%%%
%%%%%%%%%%%%%%%%%%%%%%%%%%%%%%%%%%%%%%%%%%%%%%%%%%%%%%%%%%%%%%%%%%%%%%%%%%%%%%%%%%%%%%%

%\newpage

In conclusion we state a series of questions, which can be helpful for better understanding of the structure of the twisted cnjugacy classes.

{\bf Questions.}
5. Describe the groups $G$, where $[e]_{\varphi} \leq G$ if

a) $\varphi \in \mathrm{End}\, G$ an arbitrary endomorphism;

b) $\varphi \in \mathrm{Aut}\, G$ an arbitrary automorphism;

c) $\varphi \in \mathrm{Inn}\, G$ an arbitrary inner automorophism.

As it is stated in the proposition 12, for an arbitrary group $G$ the twisted conjugacy class $[e]_{\varphi}$ is a subgroup of the group $G$ for every central automorphism $\varphi$.

6. Let $Z(G)$ be the center of the group $G$.
There are the natural maps of the group
$ \mathrm{Aut}\, G $ into a groups
$ \mathrm{Aut}\, \left( G/ Z(G)\right)$ and $ \mathrm{Aut}\, Z(G)$, for which the images of the automorphism $\varphi$ are an induced automorphism  $\overline{\varphi}$ of the factor group $ G / Z(G)$ and the restriction $\mathrm{Res} \,\varphi$
on the group $Z(G)$ respectively.
By the proposition 12 the class $[e]_{\mathrm{Res} \,\varphi}$ is always a subgroup

If $[e]_{\varphi}\leq G$, then what can we say about the factor group
$$
[e]_{\varphi}/[e]_{\mathrm{Res} \,\varphi}?
$$

Can this factor group be associated with $[\overline{e}]_{\overline{\varphi}}$,
where $\overline{e}$ the unit element of the factor group $ G/ Z(G)$?

7.  Let $N\trianglelefteq G$, $h\in G$.
We denote
$$
[e]_{h,N}=\left\{ \, [x,h] \, \| \, x\in N  \,  \right\}.
$$
It is clear that $[e]_{h,G}=[e]_{h}$.

a) Let $[e]_{h}\leq G$ for every $h\in G$.
Is it true that $[e]_{h,N} \leq N$, if $h\in G$ (if $h\in N$)?

b) How can we associate the sets $N \bigcap [e]_{h}$ and $[e]_{h,N}$
if $h\in G$ (if $h\in N$)?

8. Does there exists such a group $G$, where the set $[e]_{\varphi}$ is a subgroup for every automorphism $\varphi$ of the group $G$, but for some $\theta\in \mathrm{End}\, G$ the class $[e]_{\theta}$ is not a subgroup? Proposition 14 shows that the inverse is impossible.

%%%%%%%%%%%%%%%%%%%%%%%%%%%%%%%%%%%%%%%%%%%%%%%%%%%%%%%%%%%%%%%%%%%%%%%%%%%%%%%%%%%%%%%
%%%%%%%%%%%%%%%%%%%%%%%%%%%%%%%%%%%%%%%%%%%%%%%%%%%%%%%%%%%%%%%%%%%%%%%%%%%%%%%%%%%%%%%
%%%%%%%%%%%%%%%%%%%%%%%%%%%%%%%%%%%%%%%%%%%%%%%%%%%%%%%%%%%%%%%%%%%%%%%%%%%%%%%%%%%%%%%
%%%%%%%%%%%%%%%%%%%%%%%%%%%%%%%%%%%%%%%%%%%%%%%%%%%%%%%%%%%%%%%%%%%%%%%%%%%%%%%%%%%%%%%
%%%%%%%%%%%%%%%%%%%%%%%%%%%%%%%%%%%%%%%%%%%%%%%%%%%%%%%%%%%%%%%%%%%%%%%%%%%%%%%%%%%%%%%
%%%%%%%%%%%%%%%%%%%%%%%%%%%%%%%%%%%%%%%%%%%%%%%%%%%%%%%%%%%%%%%%%%%%%%%%%%%%%%%%%%%%%%%
%%%%%%%%%%%%%%%%%%%%%%%%%%%%%%%%%%%%%%%%%%%%%%%%%%%%%%%%%%%%%%%%%%%%%%%%%%%%%%%%%%%%%%%

%\newpage

\vspace{1cm}

%\begin{center}
%{\bf \large ЛИТЕРАТУРА }
%\end{center}

Bardakov Valery Georgievich

630090

Novosibirsk,

4 Acad. Koptyug avenue,

Sobolev institute of mathematics

2 Pirogov street,

Novosibirsk state university

E-mail: bardakov@math.nsc.ru

\vspace{1cm}

Nasybullov Timur Rinatovich

630090

Novosibirsk,

4 Acad. Koptyug avenue,

Sobolev institute of mathematics

2 Pirogov street,

Novosibirsk state university

E-mail: timur.nasybullov@mail.ru

\vspace{1cm}

Neshadim Mikhail Vladimirovich

630090

Novosibirsk,

4 Acad. Koptyug avenue,

Sobolev institute of mathematics

2 Pirogov street,

Novosibirsk state university

E-mail: neshch@math.nsc.ru


\begin{thebibliography}{99}
\bibitem{F}
 Fel'shtyn A.
Reidemeister number of any automorphism of
a Gromov hyperbolic group is infinite.
ArXiv:0101010.
\end{thebibliography}
\end{document}